# POSITIVE DYNAMIC PROGRAMMING

## A CRITIQUE AAQIB PEERZADA

**Problem Statement.** In the article "Positive Dynamic Programming", David Blackwell tries to answer the question concerning the existence of $(p, \epsilon)$-optimal stationary strategies for a positive dynamic programming problem. The principal results obtained in the paper are indicative of the existence of weakly $(p, \epsilon)$-optimal stationary plans where the optimal return need not be Borel measurable. More specifically, the main theorem [Theorem 1 in [1]], establishes the condition for weakly $(p, \epsilon)$-optimal stationary plan when the optimal return is bounded. Blackwell also gives an example that demonstrates that no $(p, \epsilon)$-optimal plan exists showing that the results obtained in his earlier work on discounted dynamic programming [2], cannot be generally extended to the positive case. Theorem 2 in [1] defines the upper bound on the income under any plan, which is an easy consequence of a theorem on the existence of upper bounds presented in [2].

**Introduction.** A dynamic programming problem is determined by four objects: $S, A, q, r$, where $S$ is a nonempty Borel set, the set of the states of the system, $A$ is the nonempty Borel set, the set of countable actions available, $q$ is the transition function that associates with each pair $(s, a)$ a probability distribution $q(.|s, a)$ on S. The transition to the next state from the current state, $s$ when an action $a$ is applied is determined by $q(.|s, a)$. $r$ is a bounded nonnegative Borel measurable function on $S \times A$, which is the immediate return. The transition of the system to a new state $s_1$ from the current state, $s$ when the action $a$ is taken gives $r(s, a, s_1)$ as the income. The process is repeated in the next state, $s_1$ and the goal is to maximize the total expected return over the infinite horizon.

A plan $\pi = (\pi_1, \pi_2, \dots )$ is a sequence where each $\pi_n$ dictates the action to take on the $n^{th}$ day based on the history of the system, $\mathcal{H}_n = (s_1, a_1, \dots, s_{n-1}, a_{n-1}, s_n)$, by associating with each $\mathcal{H}_n$ (Borel measurably) a probability distribution $\pi_n(.|\mathcal{H}_n)$ on the Borel subsets of $A$. Certain plans are of special interest. A *semi-Markov* plan is a sequence $f_1, f_2, \dots$ where each $f_n$ is a Borel measurable map from $S \times S$ into $A$, i.e. each $f_n$ is a mapping from the set of states to the set of actions, $f_n: S \to A$ and $f_n(s_1, s_n)$ is the action we take on the $n^{th}$ day if the system starts in state $s_1$ and on the $n^{th}$ day the system is in state, $s_n$. A *Markov plan* is a sequence $f_1, f_2, \dots$ where each $f_n$ is a Borel measurable map from $S \times S$ into $A$, where $f_n(s_n)$ is the action we choose on the $n^{th}$ day if the state of the system is $s_n$ on the $n^{th}$ day. A *stationary plan* is a Markov plan where $f_n = f, \forall n$ for some Borel measurable map $f: S \to A$. Such plans are denoted by $f^{(\infty)}$, and $f^{(\infty)}$ plans are stationary.

A plan $\pi$ associates with each initial state, $s$ an expected $n^{th}$ period return $r_n(\pi)(s)$ and the total expected discounted return.

$$I_\beta(\pi)(s) = \sum_{1}^{\infty} \beta^{n-1} r_n(\pi)(s) \tag{1}$$

where $\beta$ is a fixed discount factor and $0 \leq \beta < 1$. The positive dynamic programming case is the one with $r \geq 0, \beta = 1$ and the problem is of finding the plan $\pi$ that will maximize the total expected return, $I(\pi) = \sum r_n(\pi)$. The problem of positive dynamic programming was studied by David Blackwell [1], which built on the earlier results obtained in [2]. More specifically in [2], the principal results obtained were the following.

Result (I). For any probability distribution, $p$ defined on the set of states $S$, and for any $\in > 0$, there exists a stationary plan, $f^{(\infty)}$ which is $(p, \in)$ −optimal such that

$$p\{I_\beta(f^{(\infty)}) > I_\beta(\pi) - \epsilon\} = 1 \tag{2}$$

Results (II). Let $u$ be a bounded function on $S \times A$, that satisfies

$$u(s) \geq \int [r(s, a, \cdot) + \beta u(\cdot)] \, dq(\cdot \mid s, a) \tag{3}$$

Then the function $u$ is an upper bound on incomes, $I_\beta(\pi) \leq u, \forall \pi$

Result (III). If $A$ is countable, the optimal return $u_\beta^*$ is the unique bounded fixed point of the operator, $U_\beta: U \to U$, where $U$ is the set of bounded functions on $S$, defined by

$$U_\beta u(s) = \sup_a \int [r(s, a, \cdot) + \beta u(\cdot)] \, dq(\cdot \mid s, a) \tag{4}$$

**Results.** The principal results obtained in [1] are the following.

Result (i). For any probability distribution, $p$ defined on the set of states $S$, and for any $\epsilon > 0$, there *need not* exist a stationary plan, $f^{(\infty)}$ which is $(p, \epsilon)$ −optimal such that

$$p\{I(f^{(\infty)}) > I(\pi) - \epsilon\} = 1 \tag{5}$$

and the optimal return *need not* be Borel measurable. More specifically, for any probability distribution, $p$ on $S$ for which $v = \sup_\pi \int I(\pi) dp$ is finite and any $\epsilon > 0$, there is a stationary plan $f^{(\infty)}$ which is weakly $(p, \epsilon)$-optimal

$$\int I(f^\infty) \, dp > v - \epsilon \tag{6}$$

Result (ii). Any nonnegative $u$ defined on $S \times A$ that satisfies

$$u(s) \geq \int [r(s, a, \cdot) + u(\cdot)] \, dq(\cdot \,|s, a) \tag{7}$$

is an upper bound on the incomes i.e., $I(\pi) \leq u, \forall \pi$.

Result (iii). If $A$ is countable, the optimal return $u^*$ is the smallest nonnegative fixed point of the operator, $K: U \to U$, where $U$ is the set of bounded functions on $S$, defined by

$$Ku(s) = \sup_a \int [r(s, a, \cdot) + u(\cdot)] \, dq(\cdot \,|s, a) \tag{8}$$

Also, $K^n 0 \to u^*$ as $n \to \infty$.

**Example.** Consider a sequence of primary states $p(1), p(2), \ldots$ and a sequence of secondary states $s(1), s(2), \ldots$ and the terminal state, $t$. From the secondary state $s(n)$ the system moves to the state $s(n-1)$ and the income is \$1. In the secondary state $s(1)$, the system moves to the terminal state and the immediate income is \$1. Once the terminal stage is reached, no income is further received. From (ii), the income under a Markov plan $I(\pi)(s(1))$ is bounded by

$$u(s(1)) \geq r(s(1), a^1, t) + \int u(t) dq(t|s(1), a^1) \geq 1 \tag{9}$$

Thus $I(\pi)(s(1)) = 1$ and the income from the state $s(n)$ is bounded by a nonnegative function

$$u(s(n)) \geq r(s(n), a^n, s(n-1)) + \int u(s(n-1)) \, dq(s(n-1)|s(n), a^n) \geq 1 + 1 \ldots \int u(t) \, dq(t|s(1), a^1) \geq n$$

Thus $I(\pi)(s(1)) = n$. From the primary state $p(n)$ the system can either move to the next primary state $p(n+1)$ with probability 1/2 and to the terminal state, $t$ with probability 1/2. The immediate income is zero. The other choice is to move to a secondary state $s(2^n - 1)$ with immediate income $2^n - 1$. From (ii) we can write for the primary state $p(n)$

$$u(p(n)) \geq r(p(n), a^n, p(n+1)) + \int u(p(n+1)) \, dq(p(n+1)|p(n), a^n) \tag{10}$$

$$u(p(n)) \geq \int \int \ldots \int [u(p(n+k))] \, dq(p(n+k)|p(n+k-1), a^{n+k-1}) \ldots dq(p(n+1)|p(n), a^n) \tag{11}$$

For the primary state $u(p(n+k))$, application of (ii) when the system moves to the state $s(2^{n+k} - 1)$ gives

$$u(p(n+k)) \geq r(p(n+k), a^{n+k}, s(2^{n+k} - 1)) + \int u(s(2^{n+k} - 1)) dq(s(2^{n+k} - 1)|p(n+k), a^{n+k})$$

The income received when the system moves from the state $p(n+k)$ to the secondary state $s(2^{n+k} - 1)$ is $2^{n+k} - 1$. We can use this to evaluate (11) and the result is

$$u(p(n)) \geq \frac{2^{n+k} - 1}{2^k} \geq 2^n - 2^{-k} \tag{12}$$

The income from the primary state $I(\pi)(p(n))$ is nearly $2^n$ in the limiting case. Thus, the function $u: u(p(n)) = 2^n, u(s(n)) = n, u(t) = 0$ satisfies the result (ii). Any stationary policy $f^{(\infty)}$ can yield the income $I(f^{(\infty)}) = 0, \forall p(n)$ there is a primary state $p(n_0)$ from which the stationary policy moves to the secondary state $s(2^{n_0} - 1)$ so that the income $I(f^{(\infty)}) = 2^{n_0} - 1$ at $p(n_o)$. Hence for any $p$ which assigns a positive probability to every primary state and $\in < 1$, there is no $(p, \in) -$ optimal stationary plan.

Don Ornstein [3] has shown that for a certain class of positive problems with bounded optimal returns and countable state space, there is for every $\in > 0$, an $\in -$ optimal stationary plan ; $I(f^\infty) > I(\pi) - \in$. Further, Ashok Maitra [4] has shown the existence of $\in -$ optimal stationary plans when certain topological conditions are imposed on $A, q$ and $r$. More specifically, the problem that Maitra considers is the one where $A$ is a compact metric space and $r$ is bounded, nonnegative and upper semi-continuous function on $S \times A$. Under these conditions, and assuming a weak convergence $q(\cdot | s(n), a^n) \to q(\cdot | s(0), a^0)$, then for any $\in > 0$ there exists an $\in -$optimal semi- Markov plan, that is $I(\pi) \geq v^* - \in$, where $v^*$ is the optimal return $v^* = \sup_{\pi} I(\pi)$.

# ON THE EXISTENCE OF STATIONARY OPTIMAL POLICIES
## A CRITIQUE BY AAQIB PEERZADA

**Problem Statement.** The article "On the existence of stationary optimal strategies" by Don Ornstein is chiefly concerned with answering the following question. In a dynamic programming situation or a gambling situation, are strategies that consider the history of the system any better than the strategies based only on the current state or situation. The principal results obtained by Ornstein in this paper seem to suggest that in a gambling or a dynamic programming situation with an uncountably infinite state space, the family of stationary strategies are not uniformly (nearly) optimal. More specifically, the example given [Theorem A, [3]] demonstrates that no states, other than the absorbing states, can be visited more than once by any strategy. Hence, even the history-dependent strategies are no better than stationary strategies. However, in a system or situation with a countable state space, there exists a family of stationary strategies that are nearly optimal [Theorem B, [3]]. Whether this is always true, is, however, not known, even with a countable state space.

**Introduction.** Let $X$ be the states of the system and $V_x$ be the nonempty set of gambles or transition probabilities for each $x \in X$. Each gamble $v \in V_x$ is a measure on $X$ such that if the system is in the state $x \in X$, the application of $v \in V_x$ will take the system to state $y \in X$ with probability $v(y)$. A history dependent strategy chooses the gamble $v_n$ on the $n^{th}$ day as function of the history of the system $(x_1, v_1, ..., v_{n-1}, x_n)$. With a stationary strategy, however, the choice of the gamble on the $n^{th}$ day depends only on $x_n$. Starting from any state $x \in X$, suppose we want to get to a special state, $g \in X$. Then given a strategy, $s$ let $F_s(x)$ be the probability of reaching the state $g$ from $x$ when we use $s$. Let $F(x) = \sup_s F_s(x)$. Given this description of the system, the main result of the paper can be summarized as

**Result I**. Consider a system with uncountable state space $X$ with $F(x) = 1 \; \forall x \in X$. If $s$ is any stationary strategy, then for some $x_0$ in $X$, the probability of reaching state $g$, $F_s(x_0) < \frac{1}{2}$.

A related result which can be compared with Result I states:

**Result II.** If $X$ is countable, then for each $\epsilon > 0$, there exists a stationary strategy $s$ such that $F_s(x) \geq (1-\epsilon)F(x) \; \forall x \in X$.

Result I gives some information about the finite or countable set of states. Comparing the two results tells us that for countable states the stationary strategies with the property $F_s(x) \geq (1-\epsilon)F(x) \; \forall x \in X$ are nearly optimal. This limits the choice of algorithms for choosing the stationary $s$. The result II is valid for more general systems. The objective of reaching a specific state (or goal) is replaced by the income $p(x, v, x') \geq 0$ received if the system moves to state $x'$ starting from state $x$, by choosing $v \in V_x$. In this context, $F_s(x)$ is the expected value of $\sum_{i=1}^{\infty} p(x_i, v_i, x_{i+1})$. Again, let $F(x) = \sup_s F_s(x)$.

The proof of result I as given by Ornstein is based on the following reasoning. To prove I, consider an absorbing state, $b \in X$ and a set of ordinals $\alpha, \beta, \omega$ such that $\beta < \alpha < \omega$, where $\omega$ is the first uncountable ordinal. $C_\alpha$ and $C_\beta$ are the corresponding collection of points and $X = \cup_{\alpha<\omega} C_\alpha \cup \{g\} \cup \{b\}$. For each state $x \in X$, the following gamble exists, $v(g) = 1 - (\frac{1}{2^n})$, $v \in V_x$ and $v(b) = \frac{1}{2^n}$, $v \in V_x$. Two cases are possible.

1) There exists a state, $x_0$, such that $F_s(x_0) = \inf F_s(x)$ $\forall x \in \cup_{\beta<\alpha} C_\beta$. In this scenario, for each integer $n$, the following gamble is introduced: $v(x_0) = 1 - (\frac{1}{2^n})$ and $v(b) = \frac{1}{2^n}$.

2) There does not exist a state $x$ such that $F_s(x) = \inf F_s(x)$ $\forall x \in \cup_{\beta<\alpha} C_\beta$. Then let $\lim_{n\to\infty} F_s(x_n) = \inf_s F_s(x)$. For each integer, $n$ the following gamble is introduced: $v(x_n) = 1 - 2[F_s(x_n) - \inf_s F_s(x)]$ and $v(b) = 2[F_s(x_n) - \inf_s F_s(x)]$.

It is easy to see that $F(x) = 1$ $\forall x \in X$. This can be checked by realizing that $F(x) = 1, \forall x \in C_\alpha$ if we assume that $F(x) = 1, \forall x \in C_\beta$. This is true since $\beta < \alpha$. Further for any stationary strategy, $t$ that is defined on $\cup_{\beta\leq\alpha} C_\beta$, there exists a state, $x^t \in C_\alpha$ such that $F_t(x^t) < \inf F_t(y)$, where $\inf$ is taken over all $y \in \cup_{\beta<\alpha} C_\beta$. This implies that for each ordinal $\alpha$, $\inf F_s(x) < \inf F_s(y)$, $\forall x \in \cup_{\beta\leq\alpha} C_\beta$ and $y \in \cup_{\beta<\alpha} C_\beta$. This can be equivalently stated as

$$\inf F_s(y) < \frac{1}{2}, \forall y \in \cup_{\beta<\alpha} C_\beta \qquad (1)$$

The result II leads to more interesting observations about the stationary strategies when we consider systems with finite state spaces. For the result II to be valid in a more general context, Ornstein presents a third result, referred to as Theorem C in [3]. The theorem states

**Result III.** If we assume that $F(x) < \infty. \forall x$, then the result II is still true in this more general context.

The proof of the result III is quite tedious but offers many new insights. To start with we choose a state $y \in X$ and pick $\epsilon > 0, \epsilon_1 > 0$ and $\epsilon_2 > 0$ such that the following relationships are satisfied

$$\frac{1}{1+\epsilon_2} > 1 - \epsilon_1 \; ; \; 4\epsilon < \frac{\epsilon}{2} \; ; \; 8 \frac{\epsilon_1}{1} < \epsilon \quad (\frac{2}{\epsilon}) \tag{2}$$

The following conclusions can be drawn and each will lead closer to the proof of the result III.

1) For a finite set $A$, we define the outside of $A$ as another set $B = X - A$ and a strategy $s$ such that under $s$ we stop when we are outside the set $A$ i.e. in $B$. The strategy $s$ has the property given by result II i.e. $F_s(y) \geq (1 - \epsilon_1)F(y)$, where $F(y) = \sup F_s(y)$. This can be seen if we define a strategy $s_1$ such that $F_{s_1}(y) \geq (1 - \frac{1}{4}\epsilon_1)F(y)$. Let $F_s^N(x)$ be the expected amount before time, $N$, when we start from state $x$ under $s$. A particular choice of $N$ will yield $F_{s_1}^N(y) \geq (1 - \frac{1}{2}\epsilon_1)F(y)$. If $v_n(x)$ is the probability of being in state $x$ at time $n$ when starting in state $y$ under strategy $s$. Choose a finite set $A$ such that

$$\sum_{x \notin A} v_n(x)F(x) - \frac{1}{4}\epsilon_1 (\frac{1}{N}) F(y) \leq 0, \quad \forall n \leq N \tag{3}$$

When we stop outside of the set A, the stationary strategy $s_1 = s$.

2) Consider a stationary strategy $t$ such that the property $F_t(y) \geq (1 - 2\epsilon_1)F(y)$ holds when we stop in set $B$. This is true because we are considering a finite $X$. In the paper, Ornstein presents this statement in the form of a lemma.

3) For a set $E \subset A$, under $t$, such that $(1 + \epsilon_2)F_t(x) \leq F(x)$, we can write for a stationary strategy, $t'$ that stops when we are in $E$ and agrees with $t$ when we are not in $E$. In that case we can write for $t'$, $F_{t'}(y) \geq (1 - \epsilon)F(y)$. The following observation aids in understanding (3). Let $(1 + \epsilon_2)b \leq F(y)$, $a \geq (1 - \epsilon_2)F(y)$, and $F_t(y) = a + b$, where $a$ and $b$ are the expected amounts won before and after hitting the state $E$, starting in state $y$ using $t$. In that case, $F(y) \geq a + (1 + \epsilon_2)b$. This further yields

$$(a + b) \geq (1 - 2\epsilon_1)F(y) \geq (1 - 2\epsilon_1)[a + (1 + \epsilon_2)b] \geq a + (1 + \epsilon)b - 4\epsilon_1(a + b) \tag{4}$$

Using this it is clear that $b \leq \frac{1}{2}\epsilon F_t(y)$.

4) Let $V_x, x \in A - E$ include only one gamble which is designated by the stationary strategy $t$. Using the strategy $t$, we can get a new system for which the property $F^1(x) \geq (1 - \epsilon)F(x) \; \forall x \in X$ is true. If for the system we define a new stationary strategy, $s$ such that $F_s^1(y) \geq (1 - \epsilon)F(y)$, starting at the state $y$ and using strategy $s$. Now for the second state we repeat the procedure using state $y_2$ and $\frac{1}{2}\epsilon$. This process is repeated till we get result III. For the countable case, we get the limit

$$F^n(x) > [1 - \sum_{i-1}^{n} \frac{\epsilon}{2^i}] F(x) \tag{5}$$

The third part of the paper provides a discussion on as to what might happen if $p(x, v, x') < 0$. To elaborate on this, consider a system with the following transition diagram.

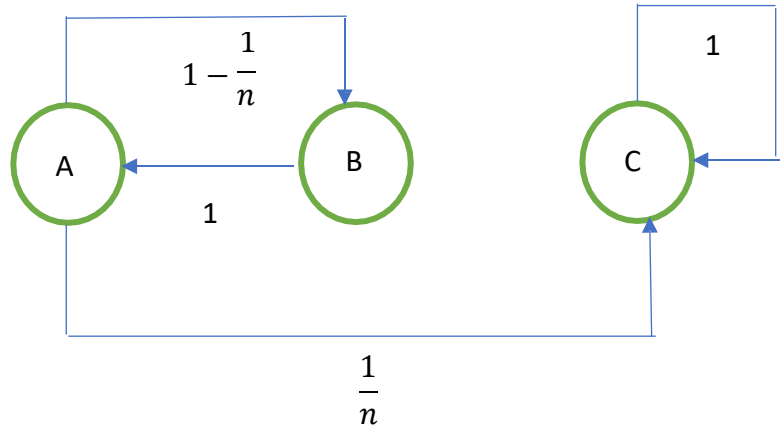

The state vector of the system, $X = [a, b, c]$. If the system starts in state $b$, it must go to state $a$ with probability 1. From state $a$, the system can move to state $c$ with probability $1/n$ or to state $b$ with probability with $1 - \frac{1}{n}$. Once the system reaches state $c$, it stays there. The income received when the system moves to state $a$ from $b$ is 0. The same is true when the system moves to state $b$ from $a$. When the system moves to state $c$, we lose a dollar and no income is received when the system is in state $c$. Any reasonable definition of $F_s$, we can write

$$F_s(a) = p(a, v(b), b) + p(a, v(c), c) = 0 - 1 = -1 \tag{6}$$

$$F(a) = \sup F_s(a) = 0 \tag{7}$$

Since $F(x)$ may be zero in case of negative $p(x, v, x')$, it doesn't make much sense to look for strategies that are good in percentage sense. Hence we aim for strategies that satisfy $F(x) - F_s(x) \geq \epsilon$. Because of this, it seems reasonable to restrict to strategies that stop with a probability 1. The conditions for optimal stationary strategies with negative $p(x, c, x')$ are

A) For systems with countable state space and bounded $\sup F_s(x)$ and $\inf F_s(x)$ (sup and inf are taken over all strategies that terminate with a probability 1), and given $\epsilon > 0$, a stationary $s$ is optimal if $F(x) - F_s(x) < \epsilon$.

B) For system with countable state space, there is a stationary $s$ satisfying $F_s(x) = F(x)$.

David Blackwell has a theorem in [2], which can be used to provide an alternate proof to the result II and hence result III as well. In [2], if the expected amount won starting on the $n^{th}$ day is discounted by a factor $\beta^n$ ($0 < \beta < 1$) where $\beta$ is the discount factor, when we start from state $x$ using $s$ and if we let $F_s^\beta(x)$ denote the expected discount won and define $F^\beta(x) = \sup_s F_s^\beta(x)$, then there exists a stationary strategy $s$ such that

$$F_s^\beta(x) \geq (1 - \epsilon) F^\beta(x) \tag{8}$$

The results in II and III follow if $\beta$ is chosen close enough to 1 so that

$$F^\beta(x) \geq (1 - \epsilon) F(x), \quad \forall x \tag{9}$$

Similar to equation (3), we can choose a gamble $v_x$ at each state $x$ such that for some strategy $s$,

$$\sum_{y \in X} F_s^\beta(y) v_x(y) + \sum_{y \in X} p(x, v_x, y) v_x(y) > \left(1 - \frac{\epsilon}{M}\right) F^\beta(x) \tag{10}$$

For a strategy $t$, if $M$ is chosen large enough, then the amount won after time $M$ will be very small and hence $t$ will be desired stationary strategy. The proof of A is very similar to that of Results III while the proof of B follows from the following arguments.

1) It is safe to assume that that each $v \in V_x$ is part of some optimal strategy.

2) If the optimal strategy exists, and $s$ is a stationary strategy such that $F_s(x) \geq (1 - \epsilon) F(x)$, then $s$ is also optimal. This follows from 1.

3) Blackwell's theorem can be used to demonstrate that $V_x$ can be replaced by a countable sub collection without changing the function $F$ or the existence of optimal strategy.

4) Each state $x$ is contained in a countable closed set and because of 2 there exists an optimal strategy for this set.

5) If an optimal strategy on a family of closed sets agrees on the intersections, then the union has an optimal stationary strategy.

Note 1: In result II, if $F(x) = \infty$ for some $(x)$, then for each integer $n$ one might expect that a stationary strategy $s$ exists such that $F_S(x) > n. \forall x$ where $F(x) = \infty$. However this is not true.

Note 2: In [2], Blackwell presents an example in which $F(x) < \infty$ but unbonded for each $(x)$ and for each stationary strategy $s$ there is an $x_0$ such that $F(x_0) - F_s(x_0) > \frac{1}{2}$.